\documentclass[twoside,11pt,leqno]{article}
\usepackage{amsfonts}

 \def\X{{\cal X}}  \def\H{{\cal H}} 

\def\i{\rm{iso}} 
\def\n{{\mathbb N}} 
\def\C{\mathbb{C}} 

\def\P{{\bf {P}}}
\def\A{\mathbb{A}} \def\e{\mathbb{B}} \def\L{\mathbb{L}} \def\R{\mathbb{R}}

\def\B{B({\cal H})} \def\b{B({\cal X})}

\newtheorem{df}{Definition}[section]
\newtheorem{thm}[df]{Theorem} \newtheorem{pro}[df]{Proposition}
\newtheorem{cor}[df]{Corollary} \newtheorem{ex}[df] {Example}
\newtheorem{rema}[df] {Remark} 
\def\sfstp{{\hskip-1em}{\bf.}{\hskip1em}}

\def\subject#1{\renewcommand{\thefootnote}{}\footnote
{AMS(MOS) subject classification (2010). Primary: {#1}}}

\def\keywords#1{\renewcommand{\thefootnote}{}\footnote
{Keywords: {#1}}}

\def\enddemo{\qed \endtrivlist} \expandafter\let\csname
enddemo*\endcsname=\enddemo

\def\qedsymbol{\ifmmode\bgroup\else$\bgroup\aftergroup$\fi
\vcenter{\hrule\hbox{\vrule
height.5em\kern.5em\vrule}\hrule}\egroup}
\def\qed{\ifmmode\else\unskip\nobreak\fi\quad\qedsymbol}

\title{\bf Isometric, symmetric and isosymmetric commuting $d$-tuples of Banach space operators}
\author{\normalsize B.P.~Duggal, I.H.~Kim}
\date{}

\begin{document}

\maketitle \thispagestyle{empty} \vskip-16pt

\subject{47A05, 47A55; Secondary47A11, 47B47.} \keywords{ Hilbert space,  Left/right multiplication operator,  $m$-left invertible, $m$-isometric and $m$-selfadjoint operators, product of operators, perturbation by nilpotents, commuting operators.  }
\footnote{The second named author was supported by Basic Science Research Program through the National Research Foundation of Korea (NRF)
 funded by the Ministry of Education (NRF-2019R1F1A1057574)}
\begin{abstract} Generalising the definition to commuting $d$-tuples of operators, a number of authors have considered structural properties of $m$-isometric, $n$-symmetric and $(m,n)$-isosymmetric commuting $d$-tuples in the recent past. This note is an attempt to take the mystique out of this extension and show how a large number of these properties follow from the more familiar arguments used to prove the single operator version of these properties.
\end{abstract}


\section {\sfstp Introduction} Let $\B$ (resp., $\b$) denote the algebra of operators, i.e. bounded linear transformations, on an infinite dimensional complex Hilbert space $\H$ into itself (resp., on an infinite dimensional comples Banach space $\X$ into itself), $\mathbb C$ denote the complex plane, $\B^d$ (resp., $\b^d$ and  ${\mathbb C}^d$) the product of $d$ copies of $\B$ (resp., $\b$ and $\mathbb C$) for some integer $d\geq 1$, $\overline{z}$ the conjugate of $z\in\mathbb C$ and ${\bf z}=(z_1,z_2,...,z_d)\in{\mathbb C}^d$. For a given polynomial $\P$ in ${\mathbb C}^d$ and a $d$-tuple ${\A}$ of commuting operators in $\B^d$, ${\A}$ is a hereditary root of $\P$ if $\P({\A})=0$. Two particular operator classes of hereditary roots which have been studied extensively are those of $m$-symmetric (also called $m$-selfadjoint in the literature) and $m$-isometric operators, where $A\in\B$ is $m$-symmetric (for some integer $m\geq 1$) if
$$\sum_{j=0}^m(-1)^j\left(\begin{array}{clcr}m\\j\end{array}\right){A^*}^{(m-j)}A^j=0$$
and $A\in\B$ is $m$-isometric if
$$\sum_{j=0}^m(-1)^j\left(\begin{array}{clcr}m\\j\end{array}\right)A^{*j}A^j=0.$$
Combining these two classes, we say $A\in\B$ is an $(m,n)$-isosymmetry (equivalently, the pair $(A^*,A)$ is $(m,n)$-isosymmetric) for some integers $m,n\geq 1$ if
\begin{eqnarray*} & &  \sum_{j=0}^m(-1)^j\left(\begin{array}{clcr}m\\j\end{array}\right)A^{*j}\left(\sum_{k=0}^n(-1)^k\left(\begin{array}{clcr}n\\k\end{array}\right)A^{*(n-k)}A^k \right)A^j\\
&=& \sum_{k=0}^n(-1)^k\left(\begin{array}{clcr}n\\k\end{array}\right)A^{*(n-k)}\left(\sum_{j=0}^m(-1)^j\left(\begin{array}{clcr}m\\j\end{array}\right)A^{*j}A^j\right) A^k= 0. 
\end{eqnarray*}
It is clear that $m$-symmetric operators  arise as solutions of $P(z)=({\overline{z}}-z)^m=0$,  $m$-isometric operators arise as  solutions of $P(z)=({\overline{z}}z-1)^m=0$ and $(m,n)$-isosymmetric operators arise as solutions of $({\overline{z}}z-1)^m({\overline{z}}-z)^n=0$. The class of $m$-symmetric operators was introduced by Helton \cite{Hel} ({\it albeit} not as operator solutions of the polynomial equation $({\overline{z}}-z)^m=0$), and the class of $m$-isometric operators was introduced by Agler \cite{A}. These classes of operators, and their variants,  have since been studied by a multitude of authors, amongst them Agler and Stankus \cite{{AS1}, {AS2}, {AS3}},  Sid Ahmed \cite {OAM}, Bayart \cite{FB}, Bermudez {\it et al} \cite{{BMN}, {BMN1}, {BMMN}}, Botelho and Jamison \cite{Fb}, Duggal \cite{{BD1}, {BD3}}, Gu \cite{{G}, {G1}} and Gu and Stankus \cite{GS}, Stankus \cite{St} and Trieu Le \cite{TL}.

\

A generalisation of the $m$-isometric property of operators $A\in\B$ to commuting $d$-tuples $\A=(A_1, \cdots, A_d)\in \B^d$, $[A_i,A_j]=A_iA_j-A_jA_i=0$ for all $1\leq i,j\leq d$, is obtained as follows \cite {GR}: $\A$ is $m$-isometric if 
$$
\sum_{j=0}^m{(-1)^j\left(\begin{array}{clcr}m\\j\end{array}\right)\sum_{|\beta|=j}{\frac{j!}{{\beta}!}}\A^{*\beta}\A^{\beta}}=0,
$$
where
\begin{eqnarray*}
& &\beta=(\beta_1, \cdots, \beta_d),  \ |\beta|=\sum_{i=1}^d \beta_i, \ {\beta}!=\Pi_{i=1}^d{\beta_i}!, \\
& &\A^{\beta}=\Pi_{i=1}^d{A_i^{\beta_i}}, \ A^{*{\beta}}=\Pi_{i=1}^d{A_i^{*{\beta_i}}};
\end{eqnarray*}
$\A$ is $n$-symmetric if 
$$
\sum_{j=0}^m{(-1)^j\left(\begin{array}{clcr}m\\j\end{array}\right)(A_1^*+ \cdots +A^*_d)^{n-j}(A_1+ \cdots +A_d)^j}=0.
$$ 
These generalisations, and certain of their variants (including $(m,n)$-isosymmetric operators), have recently been the subject matter of a number of studies, see \cite{{ACL}, {CA}, {CMN}, {GJR}} for further references.

\

This paper studies $(X,m)$-isometric, $(X,n)$-symmetric and $(X,(m,n))$-isosymmetric commuting Banach space $d$-tuples from the point of view of operators defined by elementary operators (of left and right multiplication) and shows how the arguments from the single operator case work just as well in proving a number of the structural properties of these classes of operators. The plan of the paper is as follows. In Section 2, we introduce our generalised definition of  $(X,m)$-isometric, $(X,n)$-symmetric and $(X,(m,n))$-isosymmetric commuting  $d$-tuples  in ${\b}^d$, and prove some well  known and some not so well known (possibly new) 
results on the structure of these operators. Section 3 considers perturbation by commuting nilpotent $d$-tuples, and Section 4 considers commuting products.

\section {\sfstp Definitions and introductory properties} For $A,B\in\b$, let $L_A$ and $R_B \in B(\b)$ denote respectively the operators
$$
L_A(X)=AX \ {\rm and}\ R_B(X)=XB
$$
of left multiplication by $A$ and right multiplication by $B$. A $d$-tuple $\A=(A_1, \cdots, A_d)\in \b^d$ is a commuting $d$-tuple if 
$$
[A_i,A_j]=A_iA_j-A_jA_i=0, \ {\rm all} \ 1\leq i,j\leq d.
$$ Given commuting $d$-tuples $\A =(A_1, \cdots, A_d)$ and $\e =(B_1, \cdots, B_d)$, define operators $\L_{\A}$ and $\R_{\e}$ by
$$ \L_{\A}^{\alpha}=\Pi_{i=1}^d{L^{\alpha_i}_{A_i}},  \  \R_{\e}^{\alpha}=\Pi_{i=1}^d{R_{B_i}^{\alpha_i}}
$$
where
$$
 \alpha=\left(\alpha_1, \cdots, \alpha_d\right), \ |\alpha|=\sum_{i=1}^d{\alpha_i}, \  \alpha_i\geq 0 \ {\rm for \ all} \  1\leq i\leq d.
$$
For $d$-tuples $\A$ and $\e$, and an operator $X\in\b$, let ``*"  and  ``$\times$"  denote, respectively, the multiplication operations
\begin{eqnarray*}
& & (\L_{\A} * \R_{\e})^j(X)=\left(\sum_{|\alpha|=j}{{\frac{j!}{{\alpha}!}}\L_{\A}^{\alpha}\R_{\e}^{\alpha}}\right)(X)=\left(\sum_{i=1}^d{L_{A_i}R_{B_i}}\right)^j(X)\\
& & (\rm{all \ integers} \ j\geq 0, \ {\alpha}!={\alpha_1}! \cdots {\alpha_d}!)  \ {\rm and} \\
& & \left(\L_{\A}\times\R_{\e}\right)(X)=\left(\sum_{i=1}^d L_{A_i}\right)\left(\sum_{i=1}^d R_{B_i}\right)(X). 
\end{eqnarray*}
We say that the $d$-tuples $\A$ and $\e$ commute, $[\A,\e]=0$, if 
$$
[A_i,B_j]=0 \ {\rm  for  \  all} \ 1\leq i, j\leq d.
$$ Evidently,
$$ 
[\L_{\A},\R_{\e}]=0
$$
and if $[\A,\e]=0$, then
$$
[\L_{\A},\L_{\e}]=[\R_{\A},\R_{\e}]=0.
$$
A pair $(\A,\e)$ of commuting $d$-tuples $\A$ and $\e$ is said to be  $(X,m)$-isometric, $(\A,\e)\in (X,m)$-isometric, for some positive integer $m$ and operator $X\in\b$, if
\begin{eqnarray*} 
\triangle^m_{\A,\e}(X) &=& (I-\L_{\A}  *  \R_{\e})^m(X)\\
&=&  \left( \sum_{j=0}^m(-1)^j\left(\begin{array}{clcr}n\\j\end{array}\right)\left(\L_{\A}  *  \R_{\e}\right)^j\right)(X)\\
&=&  \left(  \sum_{j=0}^m(-1)^j\left(\begin{array}{clcr}n\\j\end{array}\right)\left(\sum_{i=1}^d L_{A_i}R_{B_i}\right)^j\right)(X)\\
&=&  \sum_{j=0}^m(-1)^j\left(\begin{array}{clcr}n\\j\end{array}\right)\left(\sum_{|\alpha|=j}{\frac{j!}{{\alpha}!}}\A^{\alpha}X\e^{\alpha}\right)\\
&=& 0;
\end{eqnarray*}
$(\A,\e)$ is $(X,n)$-symmetric, for some positive integer $n$ and operator $X\in\b$, if
\begin{eqnarray*} 
\delta_{\A,\e}^n(X) &=& (\L_{\A}-\R_{\e})^n(X)\\
&=&\left( \sum_{j=0}^n(-1)^{j}\left(\begin{array}{clcr}n\\j\end{array}\right) \L_{\A}^{n-j}\times\R_{\e}^j\right)(X)\\
&=& \left(\sum_{j=0}^n(-1)^{j}\left(\begin{array}{clcr}n\\j\end{array}\right)\left(\sum_{i=1}^d L_{A_i}\right)^{n-j}\left(\sum_{i=1}^d R_{B_i}\right)^j\right)(X)\\
&=& \sum_{j=0}^m(-1)^{m-j}\left(\begin{array}{clcr}n\\j\end{array}\right)\left(\sum_{i=1}^d A_i\right)^{n-j}X\left(\sum_{i=1}^d B_i\right)^j\\
&=& 0.
\end{eqnarray*}
Commuting tuples of $(X,m)$-isometric, similarly $(X,n)$-symmetric operators, share a large number of properties with their single operator counterparts. However, there are instances where a property holds for the single operator version but fails for the $d$-tuple version. For example, whereas 
$$\ \triangle^m_{A,B}(X)=0\Longleftrightarrow \triangle^m_{A^t,B^t}(X) \ {\rm for \ all \ integers} \ t\geq 1
$$ and 
$$
\triangle^m_{A,B}(X)=0\Longleftrightarrow \triangle^m_{A^{-1},B^{-1}}(X)=0 \ {\rm for \ all \ invertible} \ A \ {\rm and} \ B,
$$
these properties fail for $d$-tuples, as the following example shows.
\begin{ex}\label{ex 00} {\em If we define operators $A_i, B_i$ ($i=1,2$) by $\A=\e= \left(\frac{1}{\sqrt{2}} I,\frac{1}{\sqrt{2}} I\right)$, then $\A=(A_1,A_2)$ and $\e=(B_1,B_2)$ are commuting, invertible $2$-tuples such that $(\A,\e)$ is $1$-isometric, i.e. $(I,1)$-isometric, but neither of $(\A^2,\e^2)$,  $\A^2=\e^2=(A^2_1,A_1A_2,A_2A_1,A_2^2)$, and $(\A^{-1},\e^{-1})$, $\A^{-1}=\e^{-1}=(\sqrt{2} I, \sqrt{2} I)$, is $m$-isometric for any $m$.}
\end{ex}
In the following we show that where a property is shared by the single operator and the $d$-tuple versions, a proof of the $d$-tuple version of the result is obtained from the argument of the single operator version of the result (if not by a transliteration of the argument, then by a simple additional argument). We remark here that a number of authors have considered $(X,m)$-isometric and $(X,n)$-symmetric Hilbert space tuples with the operator $X$ replaced by a positive operator $P$. The consideration of a general operator $X$, rather than $P\geq 0$, does not involve extra argument and does not, in general, result in loss of information. Any additional information that may result from a consideration with $P\geq 0$ is usually a result of additional hypotheses on $P$, such as injectivity, which lead to additional structure on the underlying Hilbert space. We start in the following with a couple of basic observations.
The definitions imply
$$
\triangle^t_{\A,\e}(X)=\triangle^{t-m}_{\A,\e}\left(\triangle^m_{\A,\e}(X\right),  \  \delta^t_{\A,\e}(X)=\delta^{t-n}_{\A,\e}\left(\delta^n_{\A,\e}(X)\right),
$$
\begin{eqnarray*} \triangle^{t_1}_{\A,\e}\left(\delta^{t_2}_{\A,\e}(X)\right) &=& \triangle^{t_1-m}_{\A,\e}\left[\triangle^m_{\A,\e}\left(\delta^{t_2}_{\A,\e}(X)\right)\right]\\
&=&  \triangle^{t_1-m}_{\A,\e}\left[\delta^{t_2}_{\A,\e}\left(\triangle^m_{\A,\e}(X)\right)\right]\\
 &=& \triangle^{t_1}_{\A,\e}\left[\delta^{t_2-n}_{\A,\e}\left(\delta^n_{\A,\e}(X)\right)\right]
\end{eqnarray*}
and 
\begin{eqnarray*}
 \triangle^{t_1}_{\A,\e}\left(\delta^{t_2}_{\A,\e}(X)\right) &=& \triangle^{t_1-m}_{\A,\e}\left[\triangle^m_{\A,\e}\left(\delta^{t_2-n}_{\A,\e}(\delta^n_{\A,\e}(X))\right)\right]\\
&=& \triangle^{t_1-m}_{\A,\e}\left[\delta^{t_2-n}_{\A,\e}\left(\triangle^m_{\A,\e}(\delta^n_{\A,\e}(X))\right)\right]
\end{eqnarray*}
for all integers $t_1\geq m$ and $t_2\geq n$. Thus:

\begin{pro}\label{pro00} Given commuting $d$-tuples $\A, \e$ in $\b^d$and an operator $X\in\b$,
$$
(\A,\e)\in (X,m)-{\rm isometric}\Longrightarrow (\A,\e)\in (X,t_1)-{\rm isometric \ for \ all \ integers} \  t_1\geq m;
$$
$$
(\A,\e)\in (X,n)-{\rm symmetric}\Longrightarrow (\A,\e)\in (X,t_2)-{\rm symmetric \ for \ all \ integers} \ t_2\geq n;
$$
\begin{eqnarray*}
 \{(\A,\e)&\in&  (X,m)-{\rm isometric}\}\vee \{(\A,\e)\in (X,n)-{\rm symmetric}\}\\&\Longrightarrow& (\A,\e)\in (X,(t_1,t_2))-{\rm isosymmetric}
 {\rm \ for \ all \ integers \ }t_1\geq m, t_2\geq n.
\end{eqnarray*}
\end{pro}
If $(\A,\e)\in (X,m)$-isometric, then 
\begin{eqnarray*}
\triangle^m_{\A,\e}(X)=0&\Longleftrightarrow& (I-\L_{\A}  *  R_{\e})\left(\triangle^{m-1}_{\A,\e}(X)\right)=0\\
&\Longleftrightarrow& (\L_{\A}  *  \R_{\e})\left(\triangle^{m-1}_{\A,\e}(X)\right)=\triangle^{m-1}_{\A,\e}(X)\\
&\Longrightarrow& (\L_{\A}  *  \R_{\e})^2\left(\triangle^{m-1}_{\A,\e}(X)\right)=(\L_{A}  *  R_{\e})\left(\triangle^{m-1}_{\A,\e}(X)\right)=\triangle^{m-1}_{\A,\e}(X)\\
& \cdots & \\
&\Longrightarrow& (\L_{\A}  *  \R_{\e})^t\left(\triangle^{m-1}_{\A,\e}(X)\right)=\triangle^{m-1}_{\A,\e}(X).
\end{eqnarray*}
 for all integers $t\geq 0$. Since  $\L_{A} * R_{\e}$ commutes with $\triangle^{m-1}_{\A,\e}$, we also have 
$$
\triangle^{m-1}_{\A,\e}\left((\L_{\A}  *  \R_{\e})^t(X)\right)=\triangle^{m-1}_{\A,\e}(X)
$$ 
for all integers $t\geq 0$.

Again, if $(\A,\e)\in (X,n)$-symmetric, then
\begin{eqnarray*} 
 \delta^n_{\A,\e}(X)=0 &\Longleftrightarrow& (\L_{\A}-\R_{\e})\left(\delta^{m-1}_{\A,\e}(X)\right)=0\\
&\Longleftrightarrow& \L_{\A}\left(\delta^{m-1}_{\A,\e}(X)\right)=\R_{\e}\left(\delta^{n-1}_{\A,\e}(X)\right)\\
&\cdots&\\
& \Longrightarrow& \L_{\A}^t\left(\delta^{n-1}_{\A,\e}(X)\right)=\R^t_{\e}\left(\delta^{n-1}_{\A,\e}(X)\right)
\end{eqnarray*}
for all integers $t\geq 0$. Here 
\begin{eqnarray*}
\L_{\A}\left(\delta^{n-1}_{\A,\e}(X)\right)&=&\L_{\A}\left(\sum_{j=0}^{n-1}(-1)^{j}\left(\begin{array}{clcr}n-1\\j\end{array}\right) \L_{\A}^{n-1-j}\times\R_{\e}^j\right)(X)\\
&=& \left(\sum_{j=0}^{n-1}(-1)^{j}\left(\begin{array}{clcr}n-1\\j\end{array}\right) \L_{\A}^{n-j}\times\R_{\e}^j\right)(X)
\end{eqnarray*}
and 
$$
\R_{\e}\left(\delta^{n-1}_{\A,\e}(X)\right)=\left(\sum_{j=0}^{n-1}(-1)^{j}\left(\begin{array}{clcr}n-1\\j\end{array}\right) \L_{\A}^{n-1-j}\times\R_{\e}^{j+1}\right)(X).
$$

\begin{pro}\label{pro01} Given commuting $d$-tuples $\A, \e\in\b^d$ and an operator $X\in\b$, if $(\A,\e)\in (X,m)$-isometric, then
$$
\lim_{t\rightarrow\infty}{\frac{1}{\left(\begin{array}{clcr}t\\m-1\end{array}\right)}} (\L_{\A}*\R_{\e})^t(X)=\triangle^{m-1}_{\A,\e}(X).
$$
In particular, if $(\L_{A}*\R_{\e})$ is invertible, then $\triangle^{m-1}_{\A,\e}(X)=0$.
\end{pro}
\begin{demo} The identity 
$$
(a-1)^m=a^m-\sum_{j=0}^{m-1}\left(\begin{array}{clcr}m\\j\end{array}\right)(a-1)^j
$$
applied to $\triangle^m_{\A,\e}(X)=(I-\L_{\A}*\R_{\e})^m(X)=0$ implies
$$
(\L_{\A}*\R_{\e})^m(X)=\sum_{j=0}^{m-1}\left(\begin{array}{clcr}m\\j\end{array}\right) \triangle^j_{\A,\e}(X).
$$
Observing
\begin{eqnarray*}(\L_{\A}*\R_{\e})^{m+1}(X) &=& \sum_{j=0}^{m-1}\left(\begin{array}{clcr}m\\j\end{array}\right) \triangle^{j+1}_{\A,\e}(X)+\sum_{j=0}^{m-1}\left(\begin{array}{clcr}m\\j\end{array}\right) \triangle^j_{\A,\e}(X)\\
&=& \left(\begin{array}{clcr}m\\m-1\end{array}\right) \triangle^m_{\A,\e}(X)+ \sum_{j=0}^{m-1}\left(\begin{array}{clcr}m+1\\j\end{array}\right) \triangle^j_{\A,\e}(X)\\
&=& \sum_{j=0}^{m-1}\left(\begin{array}{clcr}m+1\\j\end{array}\right) \triangle^j_{\A,\e}(X),
\end{eqnarray*}
an induction argument shows that
\begin{eqnarray*} (\L_{\A}*\R_{\e})^t(X)&=& \sum_{j=0}^{m-1}\left(\begin{array}{clcr}t\\j\end{array}\right) \triangle^j_{\A,\e}(X)\\
&=& \left(\begin{array}{clcr}t\\m-1\end{array}\right) \triangle^{m-1}_{\A,\e}(X) + \sum_{j=0}^{m-2}\left(\begin{array}{clcr}t\\j\end{array}\right) \triangle^j_{\A,\e}(X)
\end{eqnarray*}
for all integers $t\geq m$. Since $\left(\begin{array}{clcr}t\\m-1\end{array}\right)$ is of the order of $t^{m-1}$ and $\left(\begin{array}{clcr}t\\j\end{array}\right)$, $0\leq j\leq m-2$, is of the order of $t^{m-2}$ as $t\rightarrow\infty$, we have  
\begin{eqnarray*}
\triangle^{m-1}_{\A,\e}(X) &=& \lim_{t\rightarrow\infty}{\frac{1}{\left(\begin{array}{clcr}t\\m-1\end{array}\right)}}\left[(\L_{\A}*\R_{\e})^t(X)- \sum_{j=0}^{m-2}\left(\begin{array}{clcr}t\\j\end{array}\right) \triangle^j_{\A,\e}(X)\right]\\
&=&  \lim_{t\rightarrow\infty}{\frac{1}{\left(\begin{array}{clcr}t\\m-1\end{array}\right)}}(\L_{\A}*\R_{\e})^t(X).
\end{eqnarray*}
As seen above $(\L_{\A}*\R_{\e})^t\left(\triangle^{m-1}_{\A,\e}(X)\right)=\triangle^{m-1}_{\A,\e}(X)$ for all integers $t\geq 0$. Hence, if $\L_{\A}*\R_{\e}$ is invertible, then $\triangle^{m-1}_{\A,\e}(X)=(\L_{\A}*\R_{\e})^{-t}\left(\triangle^{m-1}_{\A,\e}(X)\right)$ for all integers $t\geq 0$. Consequently, if $\L_{\A}*\R_{\e}$ is invertible, then
$$
X =  \left(\begin{array}{clcr}t\\m-1\end{array}\right) \triangle^{m-1}_{\A,\e}(X) + \sum_{j=0}^{m-2}\left(\begin{array}{clcr}t\\j\end{array}\right) (\L_{\A}*\R_{\e})^{-t}\left(\triangle^j_{\A,\e}(X)\right).
$$
Since
$$
\lim_{t\rightarrow\infty}{\frac{1}{\left(\begin{array}{clcr}t\\m-1\end{array}\right)}}\left[\sum_{j=0}^{m-2}\left(\begin{array}{clcr}t\\j\end{array}\right)(\L_{\A}*\R_{\e})^{-t}\left( \triangle^j_{\A,\e}(X)\right) - X\right]=0,
$$
we have $\triangle^{m-1}_{\A,\e}(X)=0$.
\end{demo}
The case $\e=(B_1, \cdots ,B_d)\in\B^d$, $\A=\e^*$, $X=I$ and $m=2$ of Proposition \ref{pro01} is of some interest: if $(\A^*,\A)$ is $(I,2)$-isometric and $0\notin\sigma\left({\sum_{i=1}^d|A_i|^2}\right)$, then $\A$ is a spherical isometry (i.e., $(\A^*,\A)$ is $(I,1)$-isometric). For tuples $\A\in\B^d$ such that $(\A^*,\A)$ is $(I,2)$-symmetric, we have the following analogue of the well known result that $(I,2)$-symmetric operators $A\in\B$ are self-adjoint \cite{McR}.

Recall that an operator $A\in\B$ is hyponormal if $AA^*\leq A^*A$. Hyponormal pairs $(A,B^*)$ satisfy the Putnam-Fuglede commutativity property, namely that $\delta_{A,B}(I)=0$ implies $\delta_{A^*,B^*}(I)=0$. Indeed, more is true \cite{Rad}: if $\delta^n_{A,B}(I)=0$ for hyponormal $A, B^*\in\B$ and some positive integer $n$, then $\delta_{A,B}(I)=\delta_{A^*,B^*}(I)=0$.

\begin{pro}\label{pro04} If $(\A^*,\A)$ is $(I,2)$-symmetric, then $\sum_{i=1}^dA_i$ is self-adjoint.\end{pro} 

\begin{demo} For convenience, let 
$$
\sum_{i=1}^dA_i=\sum \ {\rm and} \ \sum_{i=1}^dA^*_i={\sum}^*.
$$
The hypothesis $\delta^2_{\A^*,\A}(I)=0$ then implies
$$
{\sum}^{*2}-2{\sum}^*\sum+{\sum}^2=0.
$$
Since already
$$
0\leq \delta_{\A^*,\A}(I)^*\delta_{\A^*,\A}(I)={\sum}^{*2}-{\sum}^*\sum-\sum{\sum}^*+{\sum}^2,
$$
we have (upon combining)
$$
\sum{\sum}^*\leq {\sum}^*\sum,
$$
i.e., $\sum$ is hyponormal. Evidently,
$$
\delta^2_{\A^*,\A}(I)=\sum_{j=0}^2 (-1)^{2-j}\left(\begin{array}{clcr}2\\j\end{array}\right){\sum}^{*j}{\sum}^{2-j}=0=\delta^2_{\sum,\sum^*}(I),
$$
the Putnam-Fuglede commutativity theorem applies and we conclude $\delta_{\sum,\sum^*}(I)=\delta_{\sum^*,\sum}(I)=0$. Thus $\sum={\sum}^*$.
\end{demo}

Proposition \ref{pro04} is a particular case of the following more general result, which for the case of the single operator says that an $(I,m)$-symmetric operator $T\in\B$, $m$ an even positive integer, is $(I,m-1)$-symmetric \cite{McR}. Let $T-\lambda I=T-\lambda$.

\begin{pro}\label{pro5} If $\A=(A_1, \cdots, A_d)\in \B^d$ satisfies $\delta^m_{\A^*,\A}(I)=0$ for some positive even integer $m$, then $\delta^{m-1}_{\A^*,\A}(I)=0$.
\end{pro} 

\begin{demo} The idea of the proof below is to reduce the problem to that of a single operator. For this, we start by determining the approximate point spectrum $\sigma_a(\A)$. Recall that a $d$-tuple ${\mathbb \lambda}=(\lambda_1, \cdots, \lambda_d)\in\C^d$ is in $\sigma_a(\A)$ if there exists a  sequence $\{x_n\}$ of unit vectros in $\H$ such that 
$$
\lim_{n\rightarrow\infty}{\sum_{i=1}^d\left\|(A_i-\lambda_i)x_n\right\|}=0 \ \left(\Longleftrightarrow \lim_{n\rightarrow\infty}\left\|(A_i-\lambda_i)x_n\right\|=0 \ {\rm for \ all} 1\leq i\leq d\right).
$$
Suppose ${\mathbb \lambda}\in\sigma_a(\A)$ and $\lim_{n\rightarrow\infty}\left\|(A_i-\lambda_i)x_n\right\|=0$ for all $1\leq i\leq d$. Then
$$
 \delta^m_{\A^*,\A}(I)=\sum_{j=0}^m (-1)^{j}\left(\begin{array}{clcr}m\\j\end{array}\right)\left(\sum_{i=1}^dA_i^*\right)^{m-j}\left(\sum_{i=1}^dA_i\right)^{j}=0
$$
implies
\begin{eqnarray*}
& & \lim_{n\rightarrow\infty} \sum_{j=0}^m (-1)^{j}\left(\begin{array}{clcr}m\\j\end{array}\right) \left\langle\left(\sum_{i=1}^d A_i^*\right)^{m-j}\left(\sum_{i=1}^d A_i\right)^{j}x_n, x_n\right\rangle\\
&=& \lim_{n\rightarrow\infty} \sum_{j=0}^m (-1)^{j}\left(\begin{array}{clcr}m\\j\end{array}\right) \left\langle\left(\sum_{i=1}^dA_i\right)^{j}x_n, \left(\sum_{i=1}^d A_i\right)^{m-j}x_n\right\rangle\\
&=&  \sum_{j=0}^m (-1)^{j}\left(\begin{array}{clcr}m\\j\end{array}\right) \left(\sum_{i=1}^d\lambda_i\right)^j\left(\sum_{i=1}^d{\overline{\lambda_i}}\right)^{m-j}\\
&=& \left(\sum_{i=1}^d{\overline{\lambda_i}}-\sum_{i=1}^d\lambda_i\right)^m\\
&=& 0.
\end{eqnarray*}
Hence $\sum_{i=1}^d\lambda_i$ is real for all ${\mathbb \lambda}\in\sigma_a(\A)$. The spectrum $\sigma(\A)=\sigma_a(\A)\cup\sigma_a(\A^*)$ being a compact subset of $\C$, there exists a real ${\mathbb \lambda}=(\lambda, \cdots, \lambda)\notin \sigma(\A)$ such that $\sum_{i=1}^d{(A_i-\lambda)}=\sum_{i=1}^d{A_i} -d{\lambda}$ is invertible.

\

Let $\sum_{i=1}^d{A_i}-d{\lambda}=A_{\lambda}$. Then
\begin{eqnarray*}
\delta^m_{\A^*,\A}(I) & = & (\L_{\A^*}-\R_{\A})^m(I)\\
&=& (\L_{(\A-{\mathbb \lambda})^*}-\R_{(\A-{\mathbb \lambda})})^m(I)\\
&=& \sum_{j=0}^m (-1)^{j}\left(\begin{array}{clcr}m\\j\end{array}\right)\left(\sum_{i=1}^dA_i^* -d{\lambda}\right)^{m-j}\left(\sum_{i=1}^dA_i -d{\lambda}\right)^{j}\\
&=&  \sum_{j=0}^m (-1)^{j}\left(\begin{array}{clcr}m\\j\end{array}\right) A^{*(m-j)}_{\lambda}A^j_{\lambda}\\
&=& \delta^m_{A^*_{\lambda},A_{\lambda}}(I)\\
&=& 0.
\end{eqnarray*}
The operator $A_{\lambda}$ being invertible,
\begin{eqnarray*}
\delta^m_{\A^*,\A}(I)=0 &\Longleftrightarrow& \delta^m_{\A^*_{\lambda},\A_{\lambda}}(I)=0\\
&\Longleftrightarrow& \sum_{j=0}^m (-1)^{j}\left(\begin{array}{clcr}m\\j\end{array}\right) A^{*(m-j)}_{\lambda}A^j_{\lambda}=0\\
&\Longleftrightarrow& \sum_{j=0}^m (-1)^{j}\left(\begin{array}{clcr}m\\j\end{array}\right) A^{*(-j)}_{\lambda}A^j_{\lambda}=0\\
& & \left({\rm multiply \ by} \ A^{*{-m}}_{\lambda} \ {\rm on \ the \ right}\right)\\
&\Longleftrightarrow& \triangle^m_{A^{*{-1}}_{\lambda},A_{\lambda}}(I)=0.
\end{eqnarray*}
Arguing as in the proof of \cite[Theorem 3]{DK1}, see also Proposition \ref{pro01}, this implies (recall: $m$ is even)
\begin{eqnarray*}
\triangle^{m-1}_{A^{*{-1}}_{\lambda},A_{\lambda}}(I)=0 &\Longleftrightarrow& \sum_{j=0}^{m-1} (-1)^{j}\left(\begin{array}{clcr}m-1\\j\end{array}\right) A^{*(-j)}_{\lambda}A^j_{\lambda}=0\\
&\Longleftrightarrow& \sum_{j=0}^{m-1} (-1)^{j}\left(\begin{array}{clcr}m-1\\j\end{array}\right) A^{*(m-1-j)}_{\lambda}A^j_{\lambda}=0\\
& & \left({\rm multiply \ on \ the \ left \ by} \ A^{*(m-1)}_{\lambda}\right)\\
&\Longleftrightarrow&\delta^{m-1}_{A^*_{\lambda},A_{\lambda}}(I)=0 \Longleftrightarrow \delta^{m-1}_{\A^*,\A}(I)=0.
\end{eqnarray*}
This completes the proof.
\end{demo} 

Given a sequence of operators $\{A_n\}\in\b$, we write
$$
A_n\stackrel{s}\longrightarrow A, \ A_n \ {\rm converges \ strongly \ to}  \ A,
$$
if 
$$
\lim_{n\rightarrow\infty}\left\|A_n-A\right\|=0.
$$
The $d$-tuple $\A_n=(A_{1n}, \cdots, A_{dn})$ converges strongly to $\A=(A_1, \cdots, A_d)$, $\A_n\stackrel{s}\longrightarrow \A$, if $A_{in}\stackrel{s}\longrightarrow A_i$ for all $1\leq i\leq d$. The following proposition is an analogue of a result on the norm closure of the class of $m$-isometric, similarly, $m$-symmetric, operators.
\begin{pro}\label{pro02} If $\A_n=(A_{1n}, \cdots, A_{dn})$ and $\e=(B_{1n}, \cdots, B_{dn})$ are sequences of $d$-tuples in $\b^d$ such that $A_{in}\stackrel{s}{\longrightarrow}A_i$ and $B_{in}\stackrel{s}{\longrightarrow}B_n$ for all $1\leq i\leq d$ and if either of $\triangle^{m_1}_{\A_n,\e_n}(X)$ and $\delta^{m_2}_{\A_n,\e_n}(X)$ equals $0$ for all $n$, then $\triangle^{m_1}_{\A,\e}\left(\delta^{m_2}_{\A,\e}(X)\right)=0$.
\end{pro}
\begin{demo} We start by proving that $\triangle^{m_1}_{\A_n,\e_n}(X)=0$ for all $n$ implies $\triangle^{m_1}_{\A,\e}(X)=0$ and   $\delta^{m_2}_{\A_n,\e_n}(X)=0$ for all $n$ implies  $\delta^{m_2}_{\A,\e}(X)=0$. The hypotheses $\A_n\stackrel{s}\longrightarrow \A$ and $\e_n\stackrel{s}\longrightarrow \e$ implies 
\begin{eqnarray*}
& & \lim_{n\rightarrow\infty}\|A_{in}-A_i\|=\lim_{n\rightarrow\infty}\|B_{in}-B_i\|=0\\
&\Longrightarrow& \lim_{n\rightarrow\infty}\left\|\sum_{i=1}^d(A_{in}-A_i)\right\|=\lim_{n\rightarrow\infty}\left\|\sum_{i=1}^d(B_{in}-B_i)\right\|=0
\end{eqnarray*}
for all $1\leq i\leq d$ and integers $j\geq 1$. Since
\begin{eqnarray*}
& & \left\|\triangle^{m_1}_{\A_n,\e_n}(X)-\triangle^{m_1}_{\A.\e}(X)\right\|\\
&\leq & \left\|\triangle^{m_1}_{\A_n,\e_n}(X)-\triangle^{m_1}_{\A_n.\e}(X)\right\|+\left\|\triangle^{m_1}_{\A_n,\e}(X)-\triangle^{m_1}_{\A.\e}(X)\right\|\\
&\leq& \sum_{j=0}^{m_1}\left(\begin{array}{clcr}m_1\\j\end{array}\right)\left(\left\|\sum^d_{i=1}A^j_{in}X(B^j_{in}-B^j_i)\right\|+\left\|\sum^d_{i=1}(A^j_{in}-A^j_i)XB^j_i\right\|\right)\\
&\leq& \sum_{j=0}^{m_1}\left(\begin{array}{clcr}m_1\\j\end{array}\right)\left(\sum^d_{i=1}\left\|A_{in}\right\|\|X\|\left\|B^j_{in}-B^j_i\right\|+\left\|A^j_{in}-A^j_i\right\|\|X\|\left\|B^j_i\right\|\right)
\longrightarrow 0 \\
& & {\rm as} \ n\rightarrow\infty,
\end{eqnarray*}
$\triangle^{m_1}_{\A_n,\e_n}(X)=0$ implies $\triangle^{m_1}_{\A,\e}(X)=0$. 
Considering next $\delta^{m_2}_{\A_n,\e_n}(X)$, we have

\begin{eqnarray*}
& & \left\|\delta^{m_2}_{\A_n,\e_n}(X)-\delta^{m_2}_{\A.\e}\right\|\\
&\leq & \left\|\delta^{m_2}_{\A_n,\e_n}(X)-\delta^{m_2}_{\A_n.\e}(X)\right\|+\left\|\delta^{m_2}_{\A_n,\e}(X)-\delta^{m_2}_{\A.\e}(X)\right\|\\
&\leq& \sum_{j=0}^{m_2}\left(\begin{array}{clcr}m_1\\j\end{array}\right)
\left(\left\|\left(\sum_{i=1}^d A_{in}\right)^{m_2-j}X\left[\left(\sum_{i=1}^d B_{in}\right)^j-\left(\sum_{i=1}^d B_i\right)^j\right]\right\| + \right.\\
& &\left.+ \left\|\left[\left(\sum_{i=1}^d A_{in}\right)^{m_2-j}-\left(\sum_{i=1}^d A_i\right)^{m_2-j}\right]X\left(\sum_{i=1}^d B_i\right)^j\right\|\right).
\end{eqnarray*}
Since $$
\left\|\left(\sum_{i=1}^d B_{in}\right)^j-\left(\sum_{i=1}^d B_i\right)^j\right\|\leq \left\|\sum^d_{i=1}(B_{in}-B_i)\right\|\left\|P\left(\sum_{i=1}^d B_{in},\sum_{i=1}^d B_i\right)\right\|
$$
for some polynomial $(P.,.)$, 
$$
\lim_{n\rightarrow\infty}\left\|\left(\sum_{i=1}^d B_{in}\right)^j-\left(\sum_{i=1}^d B_i\right)^j\right\|=0.
$$
Similarly,
$$
\lim_{n\rightarrow\infty}\left\|\left(\sum_{i=1}^d A_{in}\right)^{m_2-j}-\left(\sum_{i=1}^d A_i\right)^{m_2-j}\right\|=0.
$$
Hence
\begin{eqnarray*} 
 \delta^{m_2}_{\A_n,\e_n}(X)=0 \ {\rm for \ all} \ n
&\Longrightarrow& \lim_{n\rightarrow\infty}\left\|\delta^{m_2}_{\A_n,\e_n}(X)-\delta^{m_2}_{\A.\e}(X)\right\|=0\\
&\Longrightarrow& \delta^{m_2}_{\A,\e}(X)=0.
\end{eqnarray*}
Finally, since
\begin{eqnarray*}
\lim_{n\rightarrow\infty}\triangle^{m_1}_{\A_n,\e_n}\left(\delta^{m_2}_{\A_n,\e_n}(X)\right)
 &=& \lim_{n\rightarrow\infty}\triangle^{m_1}_{\A_n,\e_n}\left(\lim_{n\rightarrow\infty}\delta^{m_2}_{\A_n,\e_n}(X)\right)\\
&=& \lim_{n\rightarrow\infty}\delta^{m_2}_{\A_n,\e_n}\left(\lim_{n\rightarrow\infty}\triangle^{m_1}_{\A_n,\e_n}(X)\right),
\end{eqnarray*}
the proof is complete.
\end{demo}
Proposition \ref{pro02} is a generalisation of a number of extant results, amongst them \cite[Theorem 2]{ACL}.
\begin{rema}\label{rema00}{\em Let $\mathbb{U}=[U_{ij}]_{1\leq i,j\leq d}$ be a unitary operator in $\B^d$. Given a $d$-tuple $\mathbb{T}=(T_1, \cdots, T_d)\in \B^d$, define the $d$-tuple $\mathbb{S}=(S_1, \cdots, S_d)$ by $S_j=\sum_{i=1}^dU_{ji}T_i$; $1\leq j\leq d$. Then $\sum_{i=1}^dU^*_{ij}U_{ik}=1$ if $1\leq j=k\leq d$ and $0$ otherwise. \cite[Proposition 2.2]{GJR} claims that if $({\mathbb{T}^*},{\mathbb{T}})$ is $(A,m)$-isometric for some positive operator $A\in\B$, then $({\mathbb{S}^*},{\mathbb{S}})$ is $(A,m)$-isometric. This is false, even for single operators, as the following example shows.}
\end{rema}
\begin{ex}\label{ex00}{\em Consider operators $T=\left(\begin{array}{clcr} 1&1\\0&1\end{array}\right)$, $A=\left(\begin{array}{clcr} 0&0\\0&1
\end{array}\right)$ and $U=\left(\begin{array}{clcr} 0&1\\i&0\end{array}\right)$. Then $A\geq 0$, $U$ is unitary and $\triangle^2_{T^*,T}(A)=0$. However,
$$
S^*AS=\left(\begin{array}{clcr} 1&1\\1&1\end{array}\right), \ S^{*2}AS^2=\left(\begin{array}{clcr} 1&1-i\\1+i&2\end{array}\right)
$$
and $\triangle^2_{S^*,S}(A)\neq 0$.}
\end{ex}
We observe here that if $\triangle^m_{{\mathbb{T}},{\mathbb{T}}}(I)=0$,  and the $d$-tuple $\mathbb{S}$ and the unitary $\mathbb{U}$ are as in the remark above, then $\triangle^m_{{\mathbb{S}^*},{\mathbb{S}}}(I)=0$, as the following argument shows. We have:
\begin{eqnarray*}
\triangle^m_{{\mathbb{S}^*},{\mathbb{S}}}(I)&=&  \sum_{j=0}^{m} (-1)^j\left(\begin{array}{clcr}m\\j\end{array}\right)\left(\sum_{i=1}^d S_i^*S_i\right)^{j}\\
&=& \sum_{j=0}^{m} (-1)^j\left(\begin{array}{clcr}m\\j\end{array}\right)\left(\sum_{i=1}^d\left(\sum_{1\leq s,t\leq d}T^*_sU^*_{is}U_{it}T_t\right)\right)^{j}\\
&=&  \sum_{j=0}^{m} (-1)^j\left(\begin{array}{clcr}m\\j\end{array}\right)\left(\sum_{i=1}^d\left(\sum_{1\leq s\leq d}T^*_sU^*_{is}U_{is}T_t\right)\right)^{j}\\
&=&  \sum_{j=0}^{m} (-1)^j\left(\begin{array}{clcr}m\\j\end{array}\right)\left(\sum_{i=1}^d T^*_iT_i\right)^{j}\\
&=& \triangle^m_{{\mathbb{T}^*},{\mathbb{T}}}(I),
\end{eqnarray*}
since
\begin{eqnarray*}
 \sum_{i=1}^d\left (\sum_{1\leq s,t\leq d}T^*_sU^*_{is}U_{it}T_t\right)&=&\sum_{1\leq s,t\leq d}T^*_s\left(\sum_{i=1}^d U^*_{is}U_{it}\right)T_t\\
&=& \sum_{1\leq s\leq d}T^*_s\left(\sum_{i=1}^d U^*_{is}U_{is}\right)T_s=\sum_{1\leq s\leq d} T^*_sT_s.
\end{eqnarray*}
$(A_i,B_i)$ is $(X,m)$-isometric, even $(X,1)$-isometric, for all $1\leq i\leq d$ does not imply $(\A,\e)$ is $(X,m)$-isometric. Consider $A_i=B_i=I$ for all $1\leq i\leq d$, when it is seen that $(A_i,B_i)$ is $(X,1)$-isometric for all $X\in\b$ and $\triangle_{\A,\e}(X)=(d-1)X\neq 0$. The following proposition goes some way towards explaining this phenomenen.

\begin{pro}\label{pro03} (a) If $(A_i,B_i)$ is $(X,1)$-isometric for alll $1\leq i\leq d-1$, then $(\A,\e)$ is $(X,m)$-isometric if and only if $\left((d-2)I+L_{A_{d}}R_{B_{d}}\right)^m(X)=0$.
\

\noindent (b) If $(A_i,B_i)$ is $(X,1)$-symmetric for all $1\leq i\leq d-1$, then $(\A,\e)$ is $(X,m)$-symmetric if and only if $(A_d,B_d)$ is $(X,m)$-symmetric.
\end{pro}
\begin{demo} $(a)$  If $\triangle_{A_1,B_1}(X)=0$, then 
\begin{eqnarray*} 
\triangle^m_{\A,\e}(X) &=& \left(I- \L_{\A}*\R_{\e}\right)^m(X)=\left(I- \sum_{i=1}^d L_{A_i}R_{B_i}\right)^m(X)\\
&=&\left[\left(I- L_{A_1}R_{B_1}\right)-\left(\sum_{i=2}^d L_{A_i}R_{B_i}\right)\right]^m(X)\\
&=&  \sum_{j=0}^m (-1)^j\left(\begin{array}{clcr}m\\j\end{array}\right) \triangle^{m-j}_{A_1,B_1}\left(\sum_{i=2}^d  L_{A_i}R_{B_i}\right)^{j}(X)\\
&=&  \sum_{j=0}^m (-1)^j\left(\begin{array}{clcr}m\\j\end{array}\right) \left(\sum_{i=2}^d  L_{A_i}R_{B_i}\right)^{j}\left(\triangle^{m-j}_{A_1,B_1}(X)\right)\\
&=& 0
\end{eqnarray*}
for all $m-j\neq 0$, and if $j=m$, then $\triangle^m_{\A,\e}(X)=0$ if and only if 
$$
\triangle^m_{\A,\e}(X)=\left(\sum_{i=2}^d L_{A_i}R_{B_i}\right)^m(X)=0.
$$
Assume next that (also) $\triangle_{A_2,B_2}(X)=0$. Then $\left(\sum_{i=2}^d L_{A_i}R_{B_i}\right)^m(X)=0$ if and only if
\begin{eqnarray*}
 (-1)^m\left(\sum_{i=2}^d L_{A_i}R_{B_i}\right)^m(X)&=&\left[\left(I- L_{A_2}R_{B_2}\right)-\left(I+ \sum_{i=3}^d L_{A_i}R_{B_i}\right)\right]^m(X)\\
&=& \sum_{j=0}^m(-1)^j \left(\begin{array}{clcr}m\\j\end{array}\right) \left(I+ \sum_{i=3}^d L_{A_i}R_{B_i}\right)^{j}\left(\triangle^{m-j}_{A_2,B_2}(X)\right)\\
&=& 0.
\end{eqnarray*}
Thus  $\left(\sum_{i=2}^d L_{A_i}R_{B_i}\right)^m(X)=0$ if and only if $\left(I+ \sum_{i=3}^d L_{A_i}R_{B_i}\right)^{m}(X)=0$. Repeating the argument, we have eventually that 
$\left((d-3)I+\sum_{i=d-1}^d L_{A_i}R_{B_i}\right)^{m}(X)=0$ if and only if $\left((d-2) I+L_{A_d}R_{B_d}\right)^m(X)=0$. Conclusion:
$$
\triangle^m_{\A,\e}(X)=0\Longleftrightarrow \left((d-2) I+L_{A_d}R_{B_d}\right)^m(X)=0.
$$

\
\noindent $(b)$  If $\delta_{A_1,B_1}(X)=0$, then
\begin{eqnarray*} \delta_{\A,\e}(X) &=& (\L_{\A}-\R_{\e})^m(X)=\left(\sum_{i=1}^d(L_{A_i}-R_{B_i})\right)^m(X)\\
&=&\left[ \sum_{i=2}^d \delta_{A_i,B_i} +\delta_{A_1,B_1}\right]^m(X)\\
&=& \sum_{j=0}^m \left(\begin{array}{clcr}m\\j\end{array}\right) \left( \sum_{i=2}^d \delta_{A_i,B_i}\right)^{j}\left(\delta^{m-j}_{A_1,B_1}(X)\right)\\
&=& 0
\end{eqnarray*}
for all $m-j\neq 0$, and if $j=m$, then $\delta^m_{\A,\e}(X)=0$ if and only if $\left(\sum_{i=2}^d\delta_{A_i,B_i}\right)^m(X)=0$. Repeating the argument, we have eventually that $\delta^m_{\A,\e}(X)=0$ if and only if $\delta^m_{A_d,B_d}(X)=0$.
\end{demo}
Proposition \ref{pro03} subsumes \cite[Proposition 3]{ACL}, and proves that an analogous result holds for $(X,M)$-symmetric tuples.

\section {\sfstp Perturbation by commuting nilpotents} The single operator techniques extend to proving results on perturbation by commuting nilpotents of commuting tuples of operators satisfying an isometric or symmetic property. A  commuting $d$-tuple $\n=(N_1, \cdots, N_d)\in\b^d$ is an $n$-nilpotent for some positive integer $n$ if 
$$ 
\n^{\alpha}=\Pi^d_{i=1} N_i^{\alpha_i} =0
$$
for all $d$-tuples $\alpha=(\alpha_1, \cdots, \alpha_d)$ of non-negative integers $\alpha_i$ such that $|\alpha|=\sum_{i=1}^d \alpha_i=n$ and $\n^{\alpha}\neq 0$ for at least one $\alpha$ with $|\alpha|\leq n-1$. As usual, given  $d$-tuples $\A=(A_1, \cdots, A_d)$ and $\n=(N_1, \cdots, N_d)$, we define
$$
\A+\n=(A_1+N_1, \cdots, A_d+N_d).
$$
Recall that $[\A,\n]=0$ if and only if $[A_i,N_j]=0$ for all $1\leq i,j\leq d$.
\begin{thm}\label{thm05} Given commuting $d$-tuples $\A=(A_1, \cdots, A_d)$ and $\e=(B_1, \cdots, B_d)$  in $\b^d$ such that $\triangle^{m_1}_{\A,\e}\left(\delta^{m_2}_{\A,\e}(X)\right)=0$ for some positive integers $m_1$ and $m_2$, let
$$
\n_i=(N_{i1}, \cdots, N_{id}), \ i=1,2,
$$
 be two commuting $n_i$-nilpotent $d$-tuples such that
$$
[\A,\n_1]=[\e,\n_2]=0.
$$
Then 
$$
\triangle^{t_1}_{\A+\n_1,\e+\n_2}\left(\delta^{t_2}_{\A+\n_1,\e+\n_2}(X)\right)=0; \  t_i=m_i+n_1+n_2-2, \ i=1,2.
$$
\end{thm}
\begin{demo} The commutativity hypotheses on $\A$, $\e$, $\n_1$ and $\n_2$, taken alongwith the commutativity of the left and the right multiplication operators, imply
$$ 
\triangle^{s_1}_{\A+\n_1,\e+\n_2}\left(\delta^{s_2}_{\A+\n_1
,\e+\n_2}(X)\right)=\delta^{s_2}_{\A+\n_1,\e+\n_2}\left(\triangle^{s_1}_{\A+\n_1,\e+\n_2}(X)\right)
$$
for all positive integers $s_1$ and $s_2$. We prove the theorem in two steps. In the first step we let $\delta^{m_2}_{\A,\e}(X)=Y$; then $\triangle^{m_1}_{\A,\e}(Y)=0$, and we prove that $\triangle^{t_1}_{\A+\n_1,\e+\n_2}(Y)=0$. In the second step, we let $\triangle^{t_1}_{\A+\n_1,\e+\n_2}(Y)=Z$. Then 
$\triangle^{t_1}_{\A+\n_1,\e+\n_2}(Y)=\delta^{m_2}_{\A,\e}(Z)=0$, and we prove that $\delta^{t_2}_{\A+\n_1,\e+\n_2}(Z)=0$.

\

Considering $\triangle^{t_1}_{\A+\n_1,\e+\n_2}(Y)$, we have
\begin{eqnarray*}
 \triangle^{t_1}_{{\A+\n_1},{\e+\n_2}}
&=& (I- \L_{\A+\n_1} * \R_{\e+\n_2})^{t_1}\\
& = & \left[(I-\L_{\A} * \R_{\e})-\left((\L_{\n_1} * \R_{\e+\n_2})+(\L_{\A} * \R_{\n_2})\right)\right]^{t_1}\\
& = & \sum_{j=0}^{t_1}(-1)^j\left(\begin{array}{clcr}t_1\\j\end{array}\right)\triangle^{t_1-j}_{\A,\e}\left[\sum_{k=0}^j\left(\begin{array}{clcr}j\\k\end{array}\right)(\L_{\n_1} * \R_{\e+\n_2})^{j-k}(\L_{\A} * \R_{\n_2})^k\right]\\
& = & \sum_{j=0}^{t_1}(-1)^j\left(\begin{array}{clcr}t_1\\j\end{array}\right)\left[\sum_{k=0}^j\left(\begin{array}{clcr}j\\k\end{array}\right)(\L_{\n_1} * \R_{\e+\n_2})^{j-k}(\L_{\A} * \R_{\n_2})^k\right]\triangle^{t_1-j}_{\A,\e}.
\end{eqnarray*}
The operator $\n_i$ being $n_i$-nilpotent, $(\L_{\n_1} * \R_{\e+\n_2})^k=0$ for all $k\geq n_1$ and $(\L_{\A} * \R_{\n_2})^{j-k}=0$ for all $j-k\geq n_2$, or, $k\leq j-n_2$. Hence
$$
\triangle ^{t_1}_{\A+\n_1,\e+\n_2}(Y)=0
$$
for all $n_1\leq k\leq j-n_2$. This leaves us with the case  $n_1-1\geq k\geq j-n_2+1$. But then $j\leq n_1+n_2-2$  implies  $t_1-j\geq m_1+n_1+n_2-2-(n_1+n_2-2)= m_1$, and this, since $\triangle^{m_1}_{\A,\e}(X)=0$, forces $\triangle^{t_1-j}_{\A,\e}(Y)=0$. Conclusion: $\triangle^{t_1}_{\A+\n_1,\e+\n_2}(Y)=0$.

\

Consider now $\delta^{t_2}_{\A+\n_1,\e+\n_2}(Z)$. Since
\begin{eqnarray*}
 \delta^{t_2}_{\A+\n_1,\e+\n_2}=\left(\L_{\A+\n_1}-\R_{\e+\n_2}\right)^{t_2}
&=& \left((\L_{\A}-\R_{\e})+(\L_{\n_1}-\R_{\n_2})\right)^{t_2}\\
&=& \sum_{j=0}^{t_1}\left(\begin{array}{clcr}t_2\\j\end{array}\right)\left[\sum_{k=0}^{j}\left(\begin{array}{clcr}j\\k\end{array}\right)(-1)^{k}\L_{\n_1}^{j-k}\R_{\n_2}^k\right]\delta^{t_2-j}_{\A,\e},
\end{eqnarray*}
and since $\L_{\n_1}^{j-k}\R^k_{\n_2}=0$ for all $j-k\geq n_1$ and $k\geq n_2$, $\delta^{t_2}_{\A+\n_1,\e+\n_2}(Z)=0$ for all $n_2\leq k\leq j-n_1$. If $j-n_1+1\leq k\leq n_2-1$ (implies $j\leq n_1+n_2-2$), then $t_2-j= m_2+n_1+n_2-2-j\geq m_2$ and  $\delta^{t_2-j}_{\A,\e}(Z)=0$. Hence $\delta^{t_2}_{\A+\n_1,\e+\n_2}(Z)=0$, and the proof is complete.
\end{demo} 
Theorem \ref{thm05} subsume a number of extant results, amongst them \cite[Theorem 3.1]{GJR} and \cite[Theorem 3]{ACL}. The $d$-tuples $\A$ and $\e$ in the theorem, in the presence of suitable commutativity hypotheses, may be replaced by $d$-tuples $\A_i$ and $\e_i$; $i=1,2$. The argument of the proof of the theorem implies the following corollary.
\begin{cor}\label{cor05} Given commuting $d$-tuples $\A_i$, $\e_i$ and $\n_i$ in $\b^d$, $i=1,2$, such that 
$$
[\A_i,\n_1]=[\e_i,\n_2]=[\A_1,\A_2]=[\e_1,\e_2]=0,
$$
$$
 {\rm if} \ \triangle^{m_1}_{\A_1,\e_1}(X)=0, \ {\rm then} \ \triangle^{m_1+n_1+n_2-2}_{\A_1+\n_1,\e_1+\n_2}(X)=0,
$$
$$
{\rm if} \ \delta^{m_2}_{\A_2,\e_2}(X)=0, \ {\rm then} \ \delta^{m_2+n_1+n_2-2}_{\A_2+\n_1,\e_2+\n_2}(X)=0
$$
and
$$
{\rm if} \ \triangle^{m_1}_{\A_1,\e_1}\left(\delta_{\A_2,\e_2}^{m_2}(X)\right)=0, \ {\rm then} \ \triangle^{m_1+n_1+n_2-2}_{\A_1+\n_1,\e_1+\n_2}\left(\delta^{m_2+n_1+n_2-2}_{\A_2+\n_1,\e_2+\n_2}(X)\right)=0.
$$
\end{cor}
For the case in which $\A^*=\e=\mathbb{T}\in\B^d$ and $\n_1=\n_2=\n\in\B^d$ for some commuting $n$-nilpotent $d$-tuple $\n$ such that $[\A,\n]=[\e,\n]=0$, Theorem \ref{thm05} translates to:
\begin{cor}\label{cor050} $\triangle^{m_1}_{{\mathbb T}^*,{\mathbb T}}\left(\delta^{m_2}_{{\mathbb T}^*,{\mathbb T}}(X)\right)=0$ implies $\triangle^{m_1+2n-2}_{{\mathbb {T^*+N}},{\mathbb {T+N}}}\left(\delta^{m_2+2n-2}_{{\mathbb {T^*+N}},{\mathbb {T+N}}}(X)\right)=0$.
\end{cor} 

\section {\sfstp Isosymmetric products}  If $\A=(A_1, \cdots, A_{d_1})$ and ${\mathbb S}=(S_1, \cdots, S_{d_2})$ are two commuting $d_i$ tuples in $\b^{d_i}$, $i=1,2$, then the product ${\mathbb S}\A$ is the operator
$$
{\mathbb S}\A=(S_1A_1, \cdots, S_1A_{d_1},S_2A_2, \cdots, S_2A_{d_1}, \cdots, S_{d_2}A_1, \cdots, S_{d_2}A_{d_1}).
$$
The tuples $\A$ and ${\mathbb S}$ commute, $[\A,{\mathbb S}]=0$, if 
$$
[A_i,S_j]=0 \ {\rm for \ all} \ 1\leq i\leq d_1 \ {\rm and} \ 1\leq j\leq d_2.
$$
If $\A, \e\in \b^{d_1}$ are commuting $d_1$-tuples and ${\mathbb{S, T}}\in\b^{d_2}$ are commuting $d_2$-tuples such that 
$$
[\A,{\mathbb S}]=[\e,{\mathbb T}]=0,
$$
then ($\L_{{\mathbb S}\A}=\L_{\A{\mathbb S}}$, $\R_{{\mathbb T}\e}=\R_{\e{\mathbb T}}$ and) 
\begin{eqnarray*} 
\L_{{\mathbb S}\A} * \R_{{\mathbb T}\e}(X) &=& \sum_{j=1}^{d_2}S_j\left[\sum_{i=1}^{d_1}A_iXB_i\right]T_j\\
&=&\sum_{j=1}^{d_2}S_j\left((\L_{\A} * \R_{\e}(X)\right)T_j\\
&=& \left(\L_{{\mathbb S}} * \R_{{\mathbb T}}\right)\left((\L_{\A} * \R_{\e})(X)\right),
\end{eqnarray*}
\begin{eqnarray*}
\triangle_{{\mathbb S}\A,{\mathbb T}\e}(X) &=& (I-\L_{{\mathbb S}\A} * \R_{{\mathbb T}\e})(X)\\
&=& (I-\L_{\mathbb S}\L_{\A} * \R_{\e}\R_{\mathbb T})(X)\\
&=& \left((\L_{\mathbb S} * \R_{\mathbb T})(I-\L_{\A} * \R_{\e})+(I-\L_{\mathbb S} * \R_{\mathbb T})\right)(X)\\ 
 &=& \left((\L_{\mathbb S} * \R_{\mathbb T})\triangle_{\A,\e}+\triangle_{{\mathbb S},{\mathbb T}}\right)(X)
\end{eqnarray*}
and
$$
\left((\L_{\mathbb S} * \R_{\mathbb T})\triangle_{\A,\e}\right)^n(X)=(\L_{\mathbb S} * \R_{\mathbb T})^n\left(\triangle^n_{\A,\e}(X)\right).
$$
Again, if $\A$, $\e$, ${\mathbb S}$, ${\mathbb T}$ are the tuples above, then 
\begin{eqnarray*}
\delta_{{\mathbb S}\A,{\mathbb T}\e}(X) &=& (\L_{{\mathbb S}\A}-R_{{\mathbb T}\e})(X)\\
&=& \left[\sum_{j=1}^{d_2}L_{S_j}\left(\sum_{i=1}^{d_1} L_{A_i}\right)- \sum_{j=1}^{d_2}R_{T_j}\left(\sum_{i=1}^{d_1} R_{B_i}\right)\right](X)\\
&=& \left[\L_{\mathbb S}\times \L_{\A} - \R_{\mathbb T}\times \R_{\e}\right](X)\\
&=&\left[\L_{\mathbb S}\times (\L_{\A}-\R_{\e})+  \R_{\e}\times (\L_{\mathbb S}-\R_{\mathbb T})\right](X)\\
&=& \left[ \L_{\mathbb S}\times\delta_{\A.,\e} + \R_{\e}\times\delta_{{\mathbb S},{\mathbb T}}\right](X),
\end{eqnarray*}
and 
\begin{eqnarray*}
\left(\L_{\mathbb S}\times \delta_{\A,\e}\right)^n (X) &=& \left[\sum_{j=1}^{d_2} L_{S_j}\left(\sum_{i=1}^{d_1} L_{A_i}-R_{B_i}\right)\right]^n(X)\\
&=& \left[ \left(\sum_{j=1}^{d_2} L_{S_j}\right)^n\left(\sum_{i=1}^{d_1}(L_{A_i}-R_{B_i})\right)^n\right](X)\\
&=& (\L_{\mathbb S}^n\times\delta^n_{\A,\e})(X)=\L^n_{\mathbb S}\times\delta^n_{\A,\e}(X)
\end{eqnarray*}
and similarly
$$
(\R_{\e}\times\delta_{{\mathbb S},{\mathbb T}})^n(X)=\R_{\e}^n\times\delta^n_{{\mathbb S},{\mathbb T}}(X).
$$

It is well known, see for example \cite{{BMN}, {DM}, {DK1}, {G}, {TL}}, that if $[A_1,A_2]=[B_1,B_2]=0$ and $\triangle^{m_i}_{A_i,B_i}(X)=0$ (similarly, $\delta^{m_i}_{A_i,B_i}(X)=0$) for $i=1,2$, then $\triangle^{m_1+m_2-1}_{A_1A_2,B_1,B_2}(X)=0$ (resp., $\delta^{m_1+m_2-1}_{A_1A_2,B_1B_2}(X)=0$). Using an argument similar in spirit to the one used to prove  Theorem \ref{thm05} (see also \cite{DK1}), we prove in the following an analogous result for products $({\mathbb S}\A,{\mathbb T}{\e})$ of commuting $d$-tuples $\A, \e, {\mathbb S}$ and ${\mathbb T}$. We remark that the order $d_i$, $i=1,2$ of the $d$-tuples plays no role in the workings of our argument: there is no loss of generality in assuming $d_1=d_2=d$.
\begin{thm}\label{thm06} Let  $\A, \e, {\mathbb S}$ and ${\mathbb T}$ be commuting $d$-tuples in $\b^d$ such that 
$$
[\A,{\mathbb S}]=[\e,{mathbb S}]=[\e,{\mathbb T}]=0.
$$
If
$$
\triangle^m_{\A,\e}(\delta^n_{\A,\e}(X))=0=\triangle^r_{{\mathbb S},{\mathbb T}}(\delta^s_{{\mathbb S},{\mathbb T}}(X))
$$
and 
$$
\triangle^r_{{\mathbb S},{\mathbb T}}(\delta^n_{{\mathbb S},{\mathbb T}}(X))=0=\triangle^m_{\A,\e}(\delta^s_{\A,\e}(X))
$$
for some positive integers $m, n, r$ and  $s$, then
$$
\triangle^{t_1}_{{\mathbb S}\A,{\mathbb T}\e}(\delta^{t_2}_{{\mathbb S}\A,{\mathbb T}\e}(X))=0,
$$
where $t_1=m+r-1$ and $t_2=n+s-1$.
\end{thm}
\begin{demo} The commutativity hypothesis $[\A,{\mathbb S}]=[\e,{\mathbb T}]=0$, taken alongwith the commutativity of the left and the right multiplication operators implies
\begin{eqnarray*}
\triangle^{n_1}_{\A,\e}\left(\delta^{n_2}_{{\mathbb S},{\mathbb T}}(X)\right)
&=& \triangle_{\A,\e}^{n_1-m}\left[\triangle_{\A,\e}^m\left(\delta^{n_2-n}_{{\mathbb S},{\mathbb T}}\left(\delta^{n}_{{\mathbb S},{\mathbb T}}(X)\right)\right)\right]\\
&=& \triangle_{\A,\e}^{n_1-m}\left[\delta^{n_2-n}_{{\mathbb S},{\mathbb T}}\left(\triangle_{\A,\e}^m\left(\delta^{n}_{{\mathbb S},{\mathbb T}}(X)\right)\right)\right]
\end{eqnarray*}
for all integers $n_1\geq m$ and $n_2\geq n$. Hence
$$
\triangle^{m}_{\A,\e}\left(\delta^{n}_{{\mathbb S},{\mathbb T}}(X)\right)=0\Longrightarrow \triangle^{n_1}_{\A,\e}\left(\delta^{n_2}_{{\mathbb S},{\mathbb T}}(X)\right)=0
$$
for all integers $n_1\geq m$ and $n_2\geq n$. Similarly,

$$
\triangle^{r}_{{\mathbb S},{\mathbb T}}\left(\delta^{n}_{\A,\e}(X)\right)=0\Longrightarrow \triangle^{n_1}_{{\mathbb S},{\mathbb T}}\left(\delta^{n_2}_{\A,\e}(X)\right)=0
$$
for all integers $n_1\geq r$ and $n_2\geq n$. 

The proof below is divided into two parts. In the first part we prove
$$
\triangle^{t_1}_{{\mathbb S}\A,{\mathbb T}\e}\left(\delta^{n}_{\A,\e}(X)\right)=0, \ t_1=m+r-1,
$$
 and in the second part we prove
$$
\delta^{t_2}_{{\mathbb S}\A,{\mathbb T}\e}\left(\triangle^{t_1}_{{\mathbb S}\A,{\mathbb T}\e}(X)\right)=0, \ t_2=n+s-1.
$$
Set $\delta^n_{\A,\e}(X)=Y$. Then $\triangle^m_{\A,\e}(Y)=0$ and 
\begin{eqnarray*}
 \triangle^{t_1}_{{\mathbb S}\A,{\mathbb T}\e}(Y)&=&\left(\triangle_{{\mathbb S},{\mathbb T}}+(\L_{\mathbb S}*\R_{\mathbb T})\triangle_{\A,\e}\right)^{t_1}(Y)\\
&=& \left(\sum_{j=0}^{t_1}\left(\begin{array}{clcr}t_1\\j\end{array}\right)(\L_{\mathbb S}*\R_{\mathbb T})^{t_1-j}\triangle^{t_1-j}_{\A,\e}\triangle^j_{{\mathbb S},{\mathbb T}}\right)(Y).
\end{eqnarray*}
The commutativity hypotheses ensure
$$
\left[(\L_{\mathbb S}*\R_{\mathbb T})^{t_1-j}\triangle^{t_1-j}_{\A,\e} \ , \triangle^j_{{\mathbb S},{\mathbb T}}\right]=0;
$$ 
hence
\begin{eqnarray*}
(\L_{\mathbb S}*\R_{\mathbb T})^{t_1-j}\triangle^{t_1-j}_{\A,\e}\triangle^j_{{\mathbb S},{\mathbb T}}(Y)&=&(\L_{\mathbb S}*\R_{\mathbb T})^{t_1-j}\left(\triangle^{t_1-j}_{\A,\e}\triangle^j_{{\mathbb S},{\mathbb T}}(Y)\right)\\
&=& (\L_{\mathbb S}*\R_{\mathbb T})^{t_1-j}\left(\triangle^j_{{\mathbb S},{\mathbb T}}\triangle^{t_1-j}_{\A,\e}(Y)\right).
\end{eqnarray*}
Since
$$
\triangle_{{\mathbb S},{\mathbb T}}^j(Y)=\triangle_{{\mathbb S},{\mathbb T}}^j\left(\delta^n_{\A,\e}(X)\right)=0
$$
 for all $j\geq r$ and
$$
\triangle^{t_1-j}_{\A,\e}(Y)=\triangle^{t_1-j}_{\A,\e}\left(\delta^n_{\A,\e}(X)\right)=0
$$
for all $t_1-j\geq m$, equivalently $j\leq t_1-m=r-1$, we have 
$$
\triangle_{{\mathbb S}\A,{\mathbb T}\e}^{t_1}(\delta^n_{\A,\e}(X))=0.
$$
Now set $\delta^s_{{\mathbb S},{\mathbb T}}(X)=M$. Then $\triangle^m_{\A,\e}(M)=0$. Arguing as above, it is seen that 
$$
\triangle_{{\mathbb S}\A,{\mathbb T}\e}^{t_1}(M)=\left[\sum_{j=0}^{t_1}\left(\begin{array}{clcr}t_1\\j\end{array}\right)(\L_{\mathbb S}*\R_{\mathbb T})^{t_1-j}\triangle^{t_1-j}_{\A,\e}\triangle^j_{{\mathbb S},{\mathbb T}}\right](M).
$$
Since $\triangle^j_{{\mathbb S},{\mathbb T}}(M)=0$ for all $j\geq r$ and  $\triangle^{t_1-j}_{\A,\e}(M)=0$ for all $t_1-j\geq m$, equivalently $j\leq r-1$, we have 
$$
\triangle_{{\mathbb S}\A,{\mathbb T}\e}^{t_1}\left(\delta^s_{{\mathbb S},{\mathbb T}}(X)\right)=0.
$$
To conclude the proof, set $\triangle_{{\mathbb S}\A,{\mathbb T}\e}^{t_1}(X)=Z$. Then $\triangle_{{\mathbb S}\A,{\mathbb T}\e}^{t_1}\left(\delta^{t_2}_{{\mathbb S}\A,{\mathbb T}\e}(X)\right)=\delta^{t_2}_{{\mathbb S}\A,{\mathbb T}\e}(Z)$ and 
\begin{eqnarray*}
\delta^{t_2}_{{\mathbb S}\A,{\mathbb T}\e}(Z) &=& (\L_{{\mathbb S}\A}-\R_{{\mathbb T}\e})^{t_2}(Z)\\
&=&\left(\L_{\mathbb S}\times\delta_{\A,\e} + \R_{\mathbb T}\times\delta_{{\mathbb S},{\mathbb T}}\right)^{t_2}(Z)\\
&=& \left[\sum_{j=0}^{t_2}\left(\begin{array}{clcr}t_2\\j\end{array}\right)\left(\L_{\mathbb S}^{t_2-j}\times\delta^{t_2-j}_{\A,\e}\right)\left( \R^j_{\mathbb T}\times\delta^j_{{\mathbb S},{\mathbb T}}\right)\right](Z).
\end{eqnarray*}
Evidently,
$$
\left[\L_{\mathbb S}\times\delta_{\A,\e}, \  \R_{\mathbb S}\times\delta_{{\mathbb S},{\mathbb T}}\right]=0.
$$
Since 
$$
R_{\mathbb S}\times\delta^j_{{\mathbb S},{\mathbb T}}(Z)=R_{\mathbb S}\times\left[\delta^j_{{\mathbb S},{\mathbb T}}\left(\triangle^{t_1}_{{\mathbb S}\A,{\mathbb T}\e}(X)\right)\right]=0 
$$
for all $j\geq s$ and 
$$
\L^{t_2-j}_{\mathbb S}\times\delta^{t_2-j}_{\A,\e}(Z)=\L^{t_2-j}_{\mathbb S}\times\left(\delta^{t_2-j}_{\A,\e}\left(\triangle^{t_1}_{{\mathbb S}\A,{\mathbb T}\e}(X)\right)\right)=0
$$
for all $t_2-j\geq n$, equivalently for all $j\leq t_2-n=s-1$, we have 
\begin{eqnarray*}
\delta^{t_2}_{{\mathbb S}\A,{\mathbb T}\e}(Z) &=& \delta^{t_2}_{{\mathbb S}\A,{\mathbb T}\e}\left(\triangle^{t_1}_{{\mathbb S}\A,{\mathbb T}\e}(X)\right)\\
&=& \triangle^{t_1}_{{\mathbb S}\A,{\mathbb T}\e}\left(\delta^{t_2}_{{\mathbb S}\A,{\mathbb T}\e}(X)\right)\\
&=& 0.
\end{eqnarray*}
\end{demo}

Theorem \ref{thm06} is a generalisation of a number of currently available results on products of operators. The theorem implies, in particular the following corollaries.
\begin{cor}\label{cor06}  If $\A, \e, {\mathbb S}$ and ${\mathbb T}$ are commuting $d$-tuples in $\b^d$ such that $[\A,{\mathbb S}]=[\e,{mathbb S}]=[\e,{\mathbb T}]=0$, then

(i) $\triangle^m_{\A,\e}(X)=\triangle^n_{{\mathbb S},{\mathbb T}}(X)=0\Longrightarrow \triangle^{m+n-1}_{\A{\mathbb S},\e{\mathbb T}}(X)=0$;

(ii) $\delta^m_{\A,\e}(X)=\delta^n_{{\mathbb S},{\mathbb T}}(X)=0\Longrightarrow \delta^{m+n-1}_{\A{\mathbb S},\e{\mathbb T}}(X)=0$.
\end{cor} 
Let $C$ be a conjugation of $\H$. (Thus $C:\H\rightarrow\H$ is a conjugate linear operator such that $C^2=I$ and $\langle Cx,y\rangle=\langle Cy,x\rangle$ for all $x,y\in\H$.) The first  part of the following corollary has been proved in \cite[Theorem 4]{ACL}.
\begin{cor}\label{cor061} Let ${\mathbb S}$ and ${\mathbb T}$ be commuting $d$-tuples in $\b^d$ such that $[{\mathbb S},{\mathbb T}]=[{\mathbb S}^*,C{\mathbb T}C]=0$. Then

(i) $\triangle^m_{{\mathbb S}^*,C{\mathbb S}C}(X)=\triangle^n_{{\mathbb T}^*,C{\mathbb T}C}(X)=0\Longrightarrow \triangle^{m+n-1}_{{{\mathbb S}^*{\mathbb T}^*},C{{\mathbb S}{\mathbb T}}C}(X)=0$;

(ii)  $\delta^m_{{\mathbb S}^*,C{\mathbb S}C}(X)=\delta^n_{{\mathbb T}^*,C{\mathbb T}C}(X)=0\Longrightarrow \delta^{m+n-1}_{{\mathbb S}^*{\mathbb T}^*,C{\mathbb S}{\mathbb T}C}(X)=0$.
\end{cor}
If $\A=A\in\b$ and $\e=B\in\b$ are single operators, then the products $\A{\mathbb S}$ and $\e{\mathbb T}$ are the $d$-tuples
$$
\A{\mathbb S}=(AS_1, \cdots, AS_d)  \ {\rm and} \ \e{\mathbb T}=(BT_1, \cdots, BT_d).
$$
 If also $[A,{\mathbb S}]=[B,{\mathbb T}]=0$, then:
\begin{cor}\label{cor062} (i) $\triangle^m_{A,B}(X)=\triangle^n_{{\mathbb S},{\mathbb T}}(X)=0 \Longrightarrow \triangle^{m+n-1}_{A{\mathbb S},B{\mathbb T}}(X)=0$;

(ii) $\delta^m_{A,B}(X)=\delta^n_{{\mathbb S},{\mathbb T}}(X)=0\Longrightarrow \delta^{m+n-1}_{A{\mathbb S},B{\mathbb T}}(X)=0$.
\end{cor}
Part (i) of the corollary is a generalisation of \cite[Theorem 6]{ACL} and part (ii) of the corollary, in so far as the authors can ascertain, is new. 

\

\noindent {\bf Tensor products $\triangle^t_{\A\otimes{\mathbb S},\e\otimes{\mathbb T}}(I)$ and $\delta^t_{\A\otimes{\mathbb S},\e\otimes{\mathbb T}}(I)$}

\

Let $\X\overline{\otimes}\X$ denote the completion, endowed with a reasonable cross norm, of the algebraic tensor product of $\X$ with itself. Let $S\otimes T$ denote the tensor product of $S\in\b$ and $T\in\b$. The tensor product of the $d$-tuples $\A=(A_1, \cdots, A_d)$  and $\e=(B_1, \cdots, B_d)$ is the $d^2$-tuple
$$
\A\otimes\e=\left(A_1\otimes B_1, \cdots, A_1\otimes B_d,A_2\otimes B_1, \cdots, A_2\otimes B_d, \cdots, A_d\otimes B_1, \cdots, A_d\otimes B_d\right).
$$
Let $\A$ and $\e$ be commuting $d$-tuples such that $\triangle^m_{\A,\e}(I)=0$ (i.e., the pair $(\A,\e)$ is $m$-isometric). Then 
\begin{eqnarray*}
\triangle^m_{\A,\e}(I)=0 &\Longleftrightarrow& \triangle^m_{\A,\e}(I)\otimes I=0 \\
&\Longleftrightarrow&   \sum_{j=0}^m (-1)^j\left(\begin{array}{clcr}m\\j\end{array}\right) \left(\sum_{i=1}^d A_iB_i\right)^j\otimes I\\
&\Longleftrightarrow& \sum_{j=0}^m(-1)^j \left(\begin{array}{clcr}m\\j\end{array}\right) \left(\sum_{i=1}^d (A_i\otimes I)(B_i\otimes I)\right)^{j}=0\\
&\Longleftrightarrow& \triangle^m_{\A\otimes I, \e\otimes I}({\mathbb I})=0, \ {\mathbb I}=I\otimes I.
\end{eqnarray*}
Similarly,
$$
\triangle^m_{\A,\e}(I)=0\Longleftrightarrow \triangle^m_{I\otimes\A,I\otimes\e}({\mathbb I})=0.
$$
Considering $m$-symmetric pairs $(\A,\e)$, we have 
\begin{eqnarray*}
 \delta^m_{\A,\e}(I)=0&\Longleftrightarrow& \sum_{j=0}^m (-1)^j\left(\begin{array}{clcr}m\\j\end{array}\right)\left[ \left(\sum_{i=1}^d A_i\right)^{m-j}\left(\sum_{i=1}^d B_i\right)^j\right]\otimes I =0\\
&\Longleftrightarrow&  \sum_{j=0}^m (-1)^j\left(\begin{array}{clcr}m\\j\end{array}\right) \left(\sum_{i=1}^d A_i\otimes I\right)^{m-j}\left(\sum_{i=1}^d B_i\otimes I\right)^j=0\\
&\Longleftrightarrow& \delta^m_{\A\otimes I, \e\otimes I}({\mathbb I})=0\\
&\Longleftrightarrow& \delta^m_{I\otimes\A,I\otimes\e}({\mathbb I})=0.
\end{eqnarray*} 
The extension of Theorem \ref{thm06} to tensor products is now almost automatic.

\begin{thm}\label{thm07} Given tuples $\A, \e, {\mathbb S}$ and ${\mathbb T}$ in $\b^d$, if

(i) $\triangle^m_{\A,\e}(I)=\triangle^n_{{\mathbb S},{\mathbb T}}(I)=0$ (resp., $\delta^m_{\A,\e}(I)=\delta^n_{{\mathbb S},{\mathbb T}}(I)=0$), then $\triangle^{m+n-1}_{\A\otimes{\mathbb S},\e\otimes{\mathbb T}}({\mathbb I})=0$ (resp., $\delta^{m+n-1}_{\A\otimes{\mathbb S},\e\otimes{\mathbb T}}({\mathbb I})=0$);

(ii) $\triangle^m_{\A,\e}\left(\delta^n_{\A,\e}(I)\right)=\triangle^r_{{\mathbb S},{\mathbb T}}(I)=\delta^s_{{\mathbb S},{\mathbb T}}(I)=0$ for some positive integers $m, n, r$ and $s$, then $\triangle^{m+r-1}_{\A\otimes{\mathbb S},\e\otimes{\mathbb T}}\left(\delta^{n+s-1}_{\A\otimes{\mathbb S},\e\otimes{\mathbb T}}({\mathbb I})\right)=0$.
\end{thm}
\begin{demo} Define operators ${\bf A}$, ${\bf B}$, ${\bf S}$ and ${\bf T}$ by
$$
{\bf A}=\A\otimes I, \ {\bf B}=\e\otimes I, \ {\bf S}=I\otimes {\mathbb S}, \ {\rm and} \ {\bf T}=I\otimes{\mathbb T}.
$$
Then $[{\bf A},{\bf S}]=[{\bf B},{\bf S}]=[{\bf B},{\bf T}]=0$,  
$$
\triangle^m_{\A,\e}(I)=\triangle^n_{{\mathbb S},{\mathbb T}}(I)=0\Longleftrightarrow \triangle^m_{{\bf A},{\bf B}}({\mathbb I})=\triangle^n_{{\bf S},{\bf T}}({\mathbb I})=0
$$
and
$$
\delta^m_{\A,\e}(I)=\delta^n_{{\mathbb S},{\mathbb T}}(I)=0\Longleftrightarrow \delta^m_{{\bf A},{\bf B}}({\mathbb I})=\delta^n_{{\bf S},{\bf T}}({\mathbb I})=0.
$$
Applying Corollary \ref{cor062},
$$
\triangle^{m+n-1}_{{\bf A}{\bf S},{\bf B}{\bf T}}({\mathbb  I})=\triangle^{m+n-1}_{\A\otimes{\mathbb S},\e\otimes{\mathbb T}}({\mathbb I})=0
$$
and
$$
\delta^{m+n-1}_{{\bf A}{\bf S},{\bf B}{\bf T}}({\mathbb I})=\delta^{m+n-1}_{\A\otimes{\mathbb S},\e\otimes{\mathbb T}}({\mathbb I})=0.
$$
This proves (i).

To prove (ii), we start by observing that 
\begin{eqnarray*}
\triangle^m_{\A,\e}\left(\delta^n_{\A,\e}(I)\right) &=& \sum_{j=0}^m (-1)^j\left(\begin{array}{clcr}m\\j\end{array}\right) \sum_{k=0}^n (-1)^{k}\left(\begin{array}{clcr}n\\k\end{array}\right)\times\\
& &\times\left[\left(\sum_{i=1}^d L_{A_i}R_{B_i}\right)^{j}\left(\left(\sum_{i=1}^d L_{A_i}\right)^{n-k}\left(\sum_{i=1}^d R_{B_i}\right)^{k}\right)\right](I)\\
&=& \sum_{j=0}^m (-1)^j\left(\begin{array}{clcr}m\\j\end{array}\right) \left(\sum_{i=1}^d (-1)^{k}\left(\begin{array}{clcr}n\\k\end{array}\right) P(\A,\e)\right)
(I),
\end{eqnarray*}
where $P(\A,\e)$ is a polynomial with entries which are constant multiples of terms of type $A_1^{\alpha_1} \cdots  A_d^{\alpha_d}B_1^{\beta_1} \cdots  B_d^{\alpha_d}$ for some non-negative integers $\alpha_i, \beta_i$ ($1\leq i\leq d$). Hence
\begin{eqnarray*}  
\triangle^m_{\A,\e}\left(\delta^n_{\A,\e}(I)\right) =0 &\Longrightarrow& \triangle^m_{\A,\e}\left(\delta^n_{\A,\e}(I)\right) \otimes I=0\\
&\Longleftrightarrow&    \triangle^m_{\A\otimes I,\e\otimes I}\left(\delta^n_{\A\otimes I,\e\otimes I}({\mathbb I})\right) =0\\
&\Longleftrightarrow& \triangle^m_{{\bf A},{\bf B}}\left(\delta^n_{{\bf A},{\bf B}}({\mathbb I})\right) =0.
\end{eqnarray*}
Again
\begin{eqnarray*}
& & \triangle^r_{{\mathbb S},{\mathbb T}}(I)=0\Longleftrightarrow \triangle^r_{{\bf S},{\bf T}}({\mathbb I})=0\Longrightarrow \delta^n_{{\bf A},{\bf B}}\left(\triangle^r_{{\bf S},{\bf T}}({\mathbb I})\right)=0,\\
& & \delta^s_{{\mathbb S},{\mathbb T}}(I)=0\Longleftrightarrow \delta^s_{{\bf S},{\bf T}}({\mathbb I})=0\Longrightarrow \triangle^m_{{\bf A},{\bf B}}\left(\delta^s_{{\bf S},{\bf T}}({\mathbb I})\right)=0
\end{eqnarray*}
and
$$
\triangle^r_{{\bf S},{\bf T}} \left(\delta^s_{{\bf S},{\bf T}}({\mathbb I})\right)=0.
$$
Since ${\bf A}$, ${\bf B}$, ${\bf S}$ and ${\bf T}$ satisfy the hypotheses of Theorem \ref{thm06}, the proof of (ii) follows. 
\end{demo}


\vskip10pt \noindent\normalsize\rm B.P. Duggal,{University of Ni\v s,
Faculty of Sciences and Mathematics,
P.O. Box 224, 18000 Ni\v s, Serbia}.

\noindent\normalsize \tt e-mail:  bpduggal@yahoo.co.uk

\vskip6pt\noindent \noindent\normalsize\rm I. H. Kim, Department of
Mathematics, Incheon National University, Incheon, 22012, Korea.\\
\noindent\normalsize \tt e-mail: ihkim@inu.ac.kr

\end{document}